\documentclass{article}
\usepackage{latexsym}
\begin{document}
\title{ Factorization of Polynomials in One Variable over the Tropical Semiring}
\author{
K.H. Kim \\Mathematics Research Group\\
Alabama State University, Montgomery, AL 36101-0271, U.S.A.\\
and Fellow, Korean Academy of Science and Technology\\
email:khkim@alasu.edu\\
\and F.W. Roush\\Mathematics Research Group\\
Alabama State University, Montgomery, AL 36101-0271, U.S.A.\\
email:roushf3@hotmail.com}

\date{}
\maketitle

\thanks{2000 Mathematics Subject Classification: Primary: 16Y60; Secondary: 15A99}

\thanks{Key words and phrases: tropical algebra, tropical polynomial, Boolean polynomial,
polynomial factorization, tropical rank}

\thanks{This paper was presented at the Annual AMS meeting January 2005
in Atlanta.}

\break

\abstract{We show factorization of polynomials in one variable over the 
tropical semiring is in general NP-complete, either if all coefficients 
are finite, or if all are either $0$ or infinity (Boolean case). We give 
algorithms for the factorization problem which are not polynomial time 
in the degree, but are polynomial time for polynomials of fixed degree. 
For two-variable polynomials we derive an irreducibility criterion 
which is almost always satisfied, even for fixed degree, and is 
polynomial time in the degree. We prove there are unique least 
common multiples of tropical polynomials, but not unique greatest common 
divisors.  We show that if two polynomials in one variable 
have a common tropical factor, then their 
eliminant matrix is singular in the tropical sense. We prove the problem
of determining tropical rank is NP-hard.}

\section{Introduction}

The tropical semiring \cite{[Pin]},\cite{[Bu]}, earlier called the 
maximin algebra, is the real numbers 
together with infinity, and an additive operation $\min (x,y)$ 
and a multiplicative operation $x+y$.  
Over the past 40 years a number of 
applications for this structure have been found such as scheduling 
industrial production, hierarchical clustering, 
asymptotic approximations in physics, and nonstandard logics.  
The 2-element Boolean algebra is a 
subsemiring of the tropical semiring, and the nonnegative 
elements of the tropical semiring form the 
structure called an incline in \cite{[CKR]}. See \cite{[Gau]}, 
\cite{[Coh1]}, \cite{[Coh2]},\cite{[KR3]},
\cite{[Zim]},\cite{[Cunn]} for additional details on history 
and applications of semirings with a 
somewhat related structure.  Recently B. Sturmfels and many 
others \cite{[Stur1]}, \cite{[Stur2]}have 
discovered and exploited a new aspect to  tropical algebras, 
using them to study properties of algebraic 
varieties in terms of valuations, power series expansions, 
and Gr\"{o}bner bases.  They discovered a way 
to define algebraic varieties over a tropical semiring 
such that many properties of algebraic varieties over 
the complex numbers, such as B\'{e}zout's theorem, 
will be valid.  The paper \cite{[Stur2]} lists 5 open 
problems in the field of tropical mathematics, and the present paper deals with the second 
of these, factorization of polynomials over a tropical semiring.

This paper has a general relationship to the work of S. Gao and 
A. Lauder \cite{[Gao1]}, \cite{[Gao2]}.  The first of 
those papers gives criteria for irreducibility of multi-variable 
polynomials over the complex numbers.  These will 
also give criteria for irreducibility of multi-variable polynomials 
over the tropical semiring, since any factorization 
over a tropical semiring can be represented by a corresponding 
factorization over the complex numbers.  
The converse does not hold, since cancellation under addition to
produce zero terms does not occur in the tropical semiring.  The second of their 
papers proves NP-completeness of the problem of decomposing 
polygons (in the plane) under 
Minkowski sum and gives an algorithm for such decompositions 
which is not polynomial time in the 
strict sense but is polynomial time in the size of the integer 
coordinates, for fixed numbers of points.

Here it is likewise proved that the problem of factorization over 
a tropical semiring is NP-complete, 
but our methods and the nature of the result is quite different.  
We consider polynomials in one variable, 
and coefficients which are either 1,2, or 3, and show that the 
satisfiability problem can be expressed 
as a problem of factorization.  By factorization we do not 
mean factorization of polynomials as
functions, which is simple in the one variable case, but 
factorization as formal expansions in
which the coefficients are operated on by the rules of tropical algebra. 
We give several descriptive results on factorization and irreducibility.  
Our results apply to polynomials in several variables in a number 
of ways, such as by substituting 
different polynomials in 1 variable for each of several 
variables.  We suggest an algorithm.

We pass to the case of polynomials over a Boolean semiring (coefficients $0$ and 
$\infty$ in the tropical 
semiring) and show this factorization problem is likewise NP-complete, give irreducibility 
criteria and an algorithm also in this case.  We also discuss least common multiple, 
greatest common divisor, and eliminants.  T. Theobald \cite{[Theo]} has recently
proved other NP-completeness results in tropical geometry.

    \section{Polynomials over the tropical algebra and Minkowski sum}
The Minkowski sum of two sets $A,B$ in $n$-space is $\{ a+b |a\in A, b\in B\}$.

It is convenient to represent a polynomial
$$a_0\oplus (a_1+x)\oplus \ldots \oplus (a_n+nx)$$
over the tropical semiring geometrically, where $\oplus$ denotes 
the additive operation $\min (x,y)$ and 
the multiplicative operation is $x\otimes y=x+y$ (ordinary 
addition over the extended real numbers).  
That is, we represent the polynomial by the set of points  $(n,a_n)$ in the plane, 
together with 
all line segments from $(n,a_n)$ to $(n+1,a_{n+1})$  and all points 
above these line segments.  When we 
multiply two polynomials, the sets of points at integer x coordinates 
in this figure will combine 
by Minkowski sum, that is, in multiplying
$$a_0\oplus (a_1+x)\oplus \ldots \oplus (a_n+nx)$$
times
$$b_0\oplus (b_1+x)\oplus \ldots \oplus (b_m+mx)$$
at each degree $k$ the coefficient is $\min (a_i+b_{k-i})$. 
This gives a diagram which agrees with Minkowski sum of the diagrams of the factors 
at integer $x$-coordinates; we take the separate 
sums $a_i+b_{k-i}$ and fill in from above.   Since Minkowski sum commutes 
with convex hull, if we take the convex hulls of these figures, they will also add under 
Minkowski sum.  This implies the following which is familiar to those who 
have worked with Newton polygons.

Proposition 1. The convex hull of the diagram of a product of polynomials 
is obtained from the convex hulls of its factors by 
starting at the degree zero term, and successive drawing
 edges to the right, which as vectors equal the edges in the 
factors, taken one by one in order of increasing slope.

This follows by rational approximation and representing the tropical
polynomials by polynomials over a field in the standard fashion, which
factor according to the edges. It also follows
geometrically by comparing $\min_i (a_i+b_{k-i}), \min_i (a_i+b_{k-i+1})$. 
Each is a minimum of a difference between a sequence $a_i$ concave upward
and a sequence $-b_{k-i}$ concave downward.  By subtracting a common linear
factor $\alpha +\beta (i)$ representing a line separating two
convex sets we can arrange that this least vertical difference 
occurs at a minimum of $a_i$ and a maximum of $-b_{k-i}$.  Now consider what happens when
we shift the latter horizontally one unit.  The least difference must occur
inbetween the minimum of $a_i$ and the maximum of $-b_{k-i}$.  Thus in the
general case the minimum occurs either for the same $a_i$ and the adjacent
$b$, or vice versa.  It becomes a little more complicated when there the
minimum is not unique.

Now recursively keep track at each step of vertices 
for the two factors which add to give the current vertex in the product.  If we 
transfer any edge as a vector from a factor to the product, its right vertex will 
represent a product term entering the minimum for the product polynomial, and for the edge of 
lower slope, this will be a minimum point in the product.  This will represent
the next unused edge in one of the two factors.

Moreover terms which do not contribute to the convex hull in the factors will
not do so in the product.

Proposition 2.  A polynomial in one variable over a tropical algebra whose 
diagram lies strictly above the line segment joining its terms of highest 
and lowest degree, except at those two points, is irreducible (except for
monomial factors).

This follows from proposition 1. The following result is also known.

Proposition 3. If a polynomial has a diagram which is already convex, then it 
factors into linear factors, and this factorization can be done in polynomial time.

This suggests that the case which is most challenging may be when the convex hull of the diagram 
of the proposed product is spanned by $3$ vertices from polynomial coefficients (with an additional
infinite point), and the diagram is above the convex hull 
except at those $3$ vertices.  This determines the degrees and endpoints of the possible factors 
(by Prop.1 there can be only two factors).  

Consider the subcase of this case 
in which the first and last coefficients of the polynomial are $0$ 
(corresponding to $1$ for polynomials over a field), and their degrees are equal, which 
we call the equal slope and degree concave case.  One question is the likelihood that a 
polynomial will be irreducible, in terms of volumes 
over bounded sets of $n$-tuples of real numbers.

Proposition 4. In the equal slope and degree concave case, the 
probability is $1$ that a random polynomial 
(given the 3 vertices spanning its convex hull) is 
irreducible.  That is, the set of factorable polynomials 
of this type has strictly lower dimension 
than the set of all polynomials of this type.

Proof.  Consider two factors as above
$$a_0\oplus (a_1+x)\oplus \ldots \oplus (a_n+nx)$$
$$b_0\oplus (b_1+x)\oplus \ldots \oplus (b_m+mx)$$
where the first and last coefficients of both factors may be taken as $0$.  It will suffice to 
show that a product polynomial with coefficients $c$ which are generic, that is, not lying in 
a finite number of sets of lower dimension to be specified in this proof, cannot be factored. 
The number of coefficients $a,b$ of the factors equals the number of unknown coefficients 
$c$ of the product, and if any coefficients did not enter nontrivially into formulas for $c$, 
this 
would give a specific set of lower dimension.  Hence all $a,b$ must enter nontrivially.  We may 
also suppose all the coefficients $a,b$ except the highest and lowest degree coefficients are 
distinct, since if not this would give a set of lower dimension.
By the concavity assumption, all these other coefficients are positive.
Let $b_i$ be the largest of any of these coefficients.  Then $b_i> a_i\ge c_i$ since $c_i$  
is the infimum of a collection of terms including $a_i+0$.  But  $b_i$ 
must be directly involved
 in the minimum formula which produces some $c_j$, otherwise the product does not depend 
on all the coefficients $a,b$.  Therefore $c_j=b_i+a_{j-i}$.    But  $c_j$ is a 
minimum of terms which also include either
$a_j+0$ or $a_{j-n}>0$  depending on whether or not $j>n$.  Therefore  
some coefficient $a$ exceeds $b_i$ which is a contradiction.$\Box$

Example.  For two concave factors with unequal slope, the probability of 
being factorable is in general positive.  Consider the product of 
$$0 \oplus (a_1+x) \oplus (0+2x)$$
$$0 \oplus (b_1+x)\oplus (b_2+2x)$$
where all $a_i,b_i>0$.
For suitable inequalities the product can be 
$$0\oplus (a_1+x)\oplus (0+2x)\oplus (b_1+3x) \oplus (b_2+4x).$$
These inequalities are 
$$a_1<b_1,a_1+b_2>b_1, 2b_1>b_2$$
where the last ensures concavity of the second factor.

Products of polynomials over a tropical ring can be pictured in terms of 
minima of antidiagonals 
of a matrix.  For instance the product in the example is the antidiagonal minimum of 
$$\pmatrix{0+0  &0+b_1 &0+b_2\cr
     a_1+0   &a_1+b_1   &a_1+b_2\cr
      0+0      &0+b_1    &0+b_2}$$

\section{NP-completeness}

Lemma 5. Let  a polynomial with $c_0=c_n=c_{2n}=0$  be of the equal slope, equal degree, 
concave factorization type and $c_{2n}$  the highest degree coefficient.  Suppose all other 
$c_i\in \{ 1,2,3 \}$.  Then a factorization exists only if one exists 
of the following form: if   $c_i=1$ then 
$(a_i,b_i)=(1,1),(1,3)$ or $(3,1)$ and therefore $c_{i+n}=1$  and conversely.  
If $c_j=c_{j+n}=2$ then $(a_j,b_j)=(2,2)$.   If $c_j$  or $c_{j+n}=3$ then $(a_j,b_j)=(3,3)$.  

Proof.  By Prop.1 we may suppose the factors have $a_0=a_n=b_0=b_n=0$ and 
other coefficients are positive.   
All other coefficients $a,b$ are at least $1$, or the product would 
have additional terms below $1$.   
 If $c_j$  or  $c_{j+n} =3$ then $a_j,b_j$  must be at least $(3,3)$ or they would reduce 
the minimum in the product below $3$.  Without loss of generality we can reduce them to 
$(3,3)$.   If $c_j=c_{j+n}=2$ then $(a_j,b_j)$  are at least $(2,2)$ and when $a_j,b_j$
occur summed with other nonzero terms we have at least $3$, so we can reduce them to $(2,2)$ 
without loss of generality.   If  $c_j=1$ then one of $a_j$  or $b_j$ is $1$.  If the other
 is at least $2$, then any  value above $2$ will not affect the products, and we can 
choose it as $3$.  If the other is strictly between $1$ and $2$, then the only way it can 
affect the product coefficients nontrivially is in another product where one of $a_k,b_k$ 
is strictly between $1$ and $2$ and the other is $1$, since $c_k$ must be $1$.  But in this
case the $1+1$ will give the minimum, so the other variable can be chosen as $3$.$\Box$

Lemma 6. For any positive integer $n$, and any positive integer $N<\frac{n^{1/8}}{2}$ 
we can find 
sequences of positive numbers $x_{ij},y_{ij},z_i<n/4$ 
such that $x_{ij}+y_{ij}=z_i$ are distinct, and all 
other sums $x_{ij}+y_{rs}, x_{ij}+x_{rs},y_{ij}+y_{rs}$  
are distinct from these and from each other, and these
sums are distinct from all 
$x_{ij},y_{ij}$; here $i,j$ range from $1$ to $N$. These sequences can be found 
in polynomial time in terms of $N,n$.

Proof.  We first choose  $z_i$  in turn as distinct numbers between $n/5, n/4$.  Then in turn 
we choose the $x_{ij}$   of size at most $n/8$ which determine the $y_{ij}$ .  As we choose 
them, we avoid any equation with a previously determined number or sum.  This means 
avoiding at most $15N^8$ numbers in the choice. For instance there are at most
$N^4$ sums $x_{ij}+y_{rs}$ and $x_{ij}$ enters at most $N^2$ of them.

 But $n/8$ exceeds $16N^8$ so at each 
stage this can be done.  Choosing each variable in turn can be done in a polynomial number 
of steps, and there are a polynomial number of variables to choose.$\Box$

The degrees involved in our problem are to be polynomial in $n$.  If we fixed the degree, then 
factoring reduces to a finite number of linear inequalities which can be solved in polynomial 
time by linear programming methods (such as homotopy methods) over the rational numbers.

Definition.  The satisfiability problem is, in Boolean variables $w_i,i=1,\ldots ,n$ and 
constants $a_{ij},b_{ij},i=1,\ldots ,m,j=1,\ldots ,m$ to solve the Boolean equation
$$\prod_j \sum_i (w_i a_{ij}+w_i^c b_{ij}) =1$$
where $c$ superscript denotes complement.  It is NP-complete, in fact almost all 
other NP-complete problems are ultimately proved so by comparison with it.  
The factors are called clauses.

Theorem 7. The problem of factoring polynomials over the tropical semiring of degree $n$,
in the equal slope, equal degree concave case, with all coefficients $0,1,2,3$ is NP-complete.

Proof. It suffices to show there is a polynomial time reduction of the 
satisfiability problem to this 
factoring problem, where $m$ is of magnitude some constant positive root of $N$.  The sets of 
numbers given in Lemma 6 will specify the numbers $c$ in the factoring problem.  The degrees 
$k=x_{ij},y_{ij}$ are to be the degrees where $c_k=1$.   Twice such a degree never 
occurs elsewhere among the sums by Lemma 6, 
and in such a degree $k$ we may assume either $a_k$ or $b_k=1$ but not both, setting
$c_{2k}=c_{2k+n}=3$.  We say that 
two such degrees, that is, degrees $x_{ij}$ or $y_{ij}$,
have the same parity if it is the same 
one of $a,b$ which is $1$.  In the degrees 
$k=z_i$, we set $c_k=2$ and  $c_{k+n}=3$.   The sums of two degrees $x_{ij},y_{rs}$ can never 
exceed $n$ by Lemma 6, so $c_{k+n}=3$ can be assumed without contradiction. 

Whenever $c_k=2$ and $c_{k+n}=3$, some pair in 
degrees $x_{ij},y_{ij}$ must be added to produce $z_k$. The 
two degrees from $x_{ij},y_{ij}$ must have opposite parity, otherwise we would be adding 
$1+3,1+3$ for each combination of an $a$ and a $b$.  After some identifications the fact 
that some such pair adding to a $z_i$ has opposite parity represents some variable $w$ being 
$1$ in the corresponding clause of the satisfiability problem.  There will be one 
clause for each $z_i$.  By using the sums of $x_{ij}+y_{rs},x_{rs}+y_{ij}$, 
which are unique, and setting $c_k=2$ 
and $c_{k+n}=3$ in one of these degrees, we can force any two pairs of variables from the
$x_{ij},y_{ij}$ and $x_{rs},y_{rs}$ to have the opposite relative parity, 
without affecting other pairs.  By using 
these unique sums and $c_k=3$ and $c_{k+n}=3$ we can force any such pairs of variables to 
have the same parity.  To not require any relationship between a pair from
$x_{ij},y_{rs}$ let $c_k=2,c_{k+n}=2$.  This can be solved for by setting
$a_k=2$ without affecting other conditions. All other coefficients $c$ are to be $3$
except $c_0=c_n=c_{2n}=0$; this imposes no constraint 
since $1+1$ will never occur except when two of
$x_{ij},y_{ij}$ are added.

This allows us to identify variables in different equations 
as being the same or as being complements of each other, 
where a single variable is associated to an 
$x$ and a $y$ which add up to be a $z$; we can relate the $x$ to the $x$ and the $y$ to 
the $y$ as being either the same or complements, and make the two pairs either to have 
identical parity or unlike parity.  By identifying pairs in the same equation we can reduce 
the number of variables to any desired extent.   By these identifications we may produce 
any satisfiability problem which has at most $N/2$ variables (allowing 
for each variable to occur complemented) in each of the $N$ clauses.   
All the steps going back and forth in this 
transformation can be done in polynomial time.

If there is a factorization, then the constructions above mean that in degrees
$z_k$ there will be some pair of variables having opposite parity, so the corresponding
clause of the satisfiability problem, after we identify variables which are forced
to have the same parity, and those forced to have opposite parity to be their
complements, must have a $1$, and the satisfiability problem is solvable. 
Conversely if the satisfiability problem is solvable then we may choose a consistent
sequence of values of the parities of $x_{ij}$ and $y_{ij}$, so that whenever $k=z_i$
$c_k=2,c_{n+k}=3$, there is some pair having opposite parity.  Then for the 
$x_{ij},y_{ij}$ we choose the $a_k,b_k$ to be $(1,3)$ or $(3,1)$ according to the parities,
and in all the sums of pairs of these
degrees $c_k, c_{k+n}$ being $2,3$ or $3,3$ are accounted for, as well as all
values $c_k=1$.  In other degrees choose
the $a,b$ according to Lemma 5.  Sums of at least two of these other $a,b$ will
never cause a conflict.  This gives a tropical factorization. $\Box$

This construction also suggests that there are cases in which there are 
exponentially essentially different factorizations of tropical 
polynomials of degree $n$ into irreducible factors, since we can 
have cases of the satisfiability problem with exponentially 
many solutions, and there will be no 
refinement into common factorizations because of the location of the 
two zero coefficients in the product.

\section{Factorization of Boolean polynomials}

Here we consider factorization of Boolean polynomials
$$c_0+c_1x+\ldots +c_nx^n, c_i\in \{ 0,1\}$$
This is the same as the question of factoring polynomials over 
the tropical ring whose coefficients are $0$ and $\infty$ into 
factors of the same type.  We assume $c_0=c_n=1$ (Boolean).
Factoring Boolean polynomials is also the same as expressing 
a 1-dimensional set of numbers from $0$ to $n$ as a nontrivial Minkowski 
sum of two subsets of the same type.

Proposition 8.  If the set of positive degrees $i$ in a Boolean 
polynomial where $c_i=1$  is not a union of sets of the form 
$\{ d_1,d_2,d_1+d_2\}$ then the Boolean polynomial is 
irreducible.  In particular if the positive degrees lie 
in a set $T$ of congruence classes modulo some integer $m$, 
such that the sum of two congruence classes is outside $T$, then 
the Boolean polynomial is irreducible.  Moreover the sets of lower 
numbers $\{ d_1, d_2 \}$ suitably ordered must form some Cartesian product.

Proof.  This follows by consideration of products of 
terms $x^{d_1},x^{d_2}$ which produce the various positive 
degree terms in the polynomials. $\Box$

We do not know whether a random Boolean polynomial of degree $n$, 
for $n$ large, is more likely to factor or be irreducible.

Theorem 9.  Factorization of Boolean polynomials of degree $n$ is NP-complete.

Proof.  We reduce the factorization problem of the previous section to a 
problem of factoring two Boolean polynomials.  We give this factorization 
problem a $2$-dimensional structure by choosing a modulus $m$ less 
than the square root of $n$, and restricting degrees of the product,
which includes factor degrees, to 
those of the form $a+bm$ where $a,b<m/2$.  That 
means that in products of two of these degrees, 
the $a$'s and the $b$'s must separately add to 
produce the result, and in effect we may as well 
consider a problem of factoring a Boolean polynomial in two variables $x,y=x^m$.  
We reverse the ordering of degrees and add a constant, in passing from a
tropical polynomial of the last section to a Boolean polynomial here.  
This means an ordered pair
$(j,c_j)$  in the tropical polynomial will become a term $x^{m(0)-j}y^{m(0)-c(j)}$ 
in the two variable Boolean polynomial.  Here $m(0)$ is chosen 
large enough to avoid negatives.  We fill 
in all degrees below these down to $0$, that is, whenever a term $x^ay^b$ occurs
in a polynomial we add all terms $x^cy^d,0\le c \le a, 0 \le d\le b$.  This process
of filling in is a multiplicative homomorphism.

Now note that if the tropical polynomial factors in terms of polynomials 
represented by ordered pairs $(j,a_j)$ and $(j,b_j)$ then this polynomial 
factors in a corresponding way, where in the factors $m(0)$ is replaced by 
$m(0)/2$; if a tropical term is factored as $(a_1+m_1x)\otimes (a_2+m_2x)$ then
the Boolean term factors as 
$$(x^{m(0)/2-m_1}y^{m(0)/2-a_1})(x^{m(0)/2-m_2}y^{m(0)/2-a_2}).$$  
Where previously we took minima of coefficients of products in each degree, 
here the maximum will be taken, since we fill in below.

Conversely if there is any factorization of the Boolean polynomial, because 
the constant term is $1$, the factors must represent subsets of the terms 
producing the product, and we can fill in all degrees from below without loss 
of generality.  Then the top degree terms must multiply in the two factors like 
the products in a tropical algebra, and must produce the given product.  Therefore 
the Boolean polynomial factors nontrivially if and only if the tropical 
polynomial does. $\Box$

\section{Tropical eliminant and tropical rank}

Given two polynomials $f(x), g(x)$ of degrees $r,s$, we form the eliminant 
matrix $E$ of degree 
$r+s$ just as over the complex numbers.  Its $r+s$ rows are the vectors of coefficients of 
$f,xf,\ldots ,x^{s-1}f, g,xg,\ldots ,x^{r-1}g$.  This produces a square matrix, and over the complex numbers 
if $t$ is a common root of $f,g$ then the matrix multiplied by the column vector $1,t,\ldots ,t^{r+s-1}$ 
 is zero therefore it is singular.  This then gives a criterion for a common factor to exist which has 
many uses. In our case, missing coefficients over the tropical ring are considered infinite.

A $n\times n$ tropical matrix $(a_{ij})$ is said to be nonsingular if and only if the
permutation $\pi$ minimizing $\sum_{i=1}^1 a_{i\pi (i)}$ is unique, otherwise
singular.

This concept is isomorphic to strong regularity in max-algebra which has been studied
in recent work of P. Butkovi\v{c}.

Theorem 10. If $f(x), g(x)$ are polynomials in one variable over the tropical semiring 
having a common factor, then their tropical eliminant is singular.

Proof. First consider the case of polynomials over the tropical rationals.
By scaling we can reduce this to matrices over the nonnegative integers.  Represent
such polynomials as ordinary polynomials in $t$ whose coefficients 
are polynomials in a parameter $\lambda$.  Assume the tropical polynomials 
have a common tropical
factor.  Then we can represent them by ordinary polynomials with a common factor, for 
which otherwise coefficients are generic.

Then for each value of $\lambda$ the ordinary eliminants must be singular.  Yet 
for $\lambda$ sufficiently  close to zero, the terms of 
lowest order in $\lambda$ will be
dominant, and if there is only one of given order, it cannot cancel.  Therefore multiple
terms of this least order must occur, and the tropical eliminant matrix must be
singular.

The extension to the tropical reals follows from the case of the tropical rationals
by taking sufficiently close approximations to the real numbers involved, which
result from a homomorphism of $Q$-vector spaces $R$ to $Q$. For such an approximations
equalities will follow from equalities in the real case, and also inequalities from
inequalities in the real case.  (Alternatively we could take general real powers
of $\lambda$ in the coefficients). $\Box$

Example.  Two polynomials of degree $n$ with $a_0=a_n=0,b_0=b_n=0$ and
other coefficients positive will have singular eliminant but will
not in general have a common factor, so this condition is not sufficient (except
possibly for the convex case). 

The following is a proposed problem in tropical mathematics not 
directly related to factoring polynomials, but somewhat to the above considerations.

Definition.  The tropical rank of a matrix 
is the largest size of its nonsingular 
submatrices.

The row (column) space of a Boolean matrix is the space spanned by its 
rows (columns) under Boolean addition and multiplication by $0$. 
It is proved by Devlin, Santos, and 
Sturmfels \cite{[Dev]}
that the tropical rank of a $(0,1)$-tropical matrix with no row (column) 
consisting entirely of $1$ is the length of a longest chain
in the row space (column space--by the next result the two are equivalent) of
a complementary Boolean matrix, omitting the zero vector.  For instance the 
rank of the tropical matrix which is 0 on the main diagonal and 1 everywhere else
is its dimension.
It is sufficient to consider only chains in the row space 
which are formed by adding one Boolean basis 
vector at a time to the existing sum: if a chain has $\sum x_i > \sum y_j$ at
sum step then replace this by $\sum x_i +\sum y_j >\sum y_j$ and then add in the
vectors $x_i$ one by one when they make a difference, 
possibly making the chain longer. In the case of two-valued matrices with two
finite values, being nonsingular is equivalent to the corresponding order-reversed
Boolean $(0,1)$-matrix either having permanent $1$ or being a direct sum of a permanent
$1$ matrix with a $1\times 1$ zero matrix: to see this consider a minimizing permutation
and the two submatrices spanned by its entries with each of the higher and
lower values.

Two Boolean
matrices $A,B$ are L (R) equivalent if and only if they have the same
row (column) space; they are D-equivalent if and only if there is a Boolean
matrix $C$ such that $A$ is R-equivalent to $C$ and $C$ is L-equivalent to
$B$.  Though we will not actually use the next two results, they show the
close connection between general finite lattices and Boolean matrices.

Theorem (Markowsky \cite{[Mar]}).  There is an isomorphism from the row space of a finite
Boolean matrix as a lattice to the inversion of the column space.  Any finite
lattice is equivalent to the row space of a Boolean matrix, which is unique up
to D-equivalence.

It is possible to explicitly construct the required isomorphism in the first
statement.  Every finite lattice can be represented in terms of its additive
operation as a semilattice of sets, looking in terms of its order ideals,
and this will be the row space of a Boolean matrix.  A proof and reference is given in
\cite{[Kim]} that two Boolean matrices are D-equivalent if and only if their
row spaces are isomorphic as lattices. This result was first proved by 
Zaretskii.

In connection with what lattices can arise, we mention a result of Markowsky,
which is likely related to Theorem 23 \cite{[DevSturm]}.

Theorem (Markowsky). A finite lattice is isomorphic to the row space of an
$n\times n$ Boolean matrix if and only if it has at most 
$n$ elements in its basis and at most $n$ elements in its dual basis.

The argument goes by constructing an actual 
Boolean matrix by letting rows correspond
to basis elements, columns to dual basis elements, and letting the $(i,j)$
entry be 1 if and only if a basis element is not less than a dual basis element.
This gives an $n\times n$ Boolean matrix whose 
rows and columns have the given relationship.
But some other, possibly larger Boolean matrix represents the finite lattice.
However the two Boolean matrices 
must be D-equivalent because they reflect the same relations
between generators and dual generators.  Dependent rows 
and columns will not affect the
row space lattice, or the D-class.

Work of P. Butkovi\v{c} and F. Hevery \cite{BuH} on strong
regularity has essentially shown that determining whether a $k\times k$
tropical matrix has maximal rank can be done in polynomial time, using
an algorithm for the assignment problem, and then studying the 0 and nonzero
values of the final matrix to determine whether a permuted diagonal of zero
entries is unique.  Another case for which a polynomial time algorithm 
exists is the following.

Proposition 11.  For $k\times n, k<n$ tropical matrices whose 
entries have two distinct values, such
that the higher of these occurs in every column there is a polynomial time algorithm
$(O(n^3))$ to decide whether they have tropical rank $k$.  

Proof. As in \cite{[Dev]} this problem is equivalent to find an increasing chain of length
$k$ in the nonzero column space of the corresponding Boolean matrix.  This chain must
begin with a vector having one $1$ entry and each time some $1$ entry must be added 
outside the previous set of $1$ entries in the sequence.  This can be decided by
first listing all distinct columns with a single 1, then in turn 
all columns having a single 1 outside the set of columns previously listed.  If this
process continues for $k$ stages then it produces, by sums over the initial segments
of the sequence, an increasing chain of length $k$.  If the process ever stops with
a set $S$ of $1$ entries then this chain $C$ cannot produce
a solution, and also any other chain $C_1$ of that length starting with a column with a
single 1 must terminate, for if not, at some point a column of $C_1$ could be inserted to
add exactly one 1 outside the columns of $C$.  That is because the set of 1 entries
in $C_1$ grows by 1 each time from a single 1 entry, to all possible 1 entries.
$\Box$

This generalizes to size $k+c$ or rank $k-c$ for each constant $c$: just 
consider all sets of $k$ rows and test them as above.  The assumption on the
lower value occurring in each column can also be removed.

By a principal triangular submatrix of size $k$ in a matrix $M$, we mean a
submatrix which is conjugated by some permutation $P$ to the locations
$(i,j)|i,j\le k$ and that in these locations $PMP^T$ is 1 above
the main diagonal.

Lemma 12. The problem of finding a maximal triangular principal submatrix
in a given $(0,1)$-matrix is NP-complete.

Proof. In an undirected digraph this is equivalent to the NP-complete
complete subgraph problem.

Theorem 13. Determining whether $n\times n$ 
Boolean matrices have tropical rank at least $k$ is in general NP-complete.

Proof. We consider this in terms of the longest chain in the column space of the
Boolean matrix.
We restrict first to a special type of Boolean matrix.  Vertices (indices of the
matrix) are
divided into 2 sets and columns have 1 vertex in the first set and
2 vertices in the second set.  This is regarded as an edge coloring 
of a graph where the last 2 vertices specify the edge and the first
vertex specifies the color. Any element
minimally greater than another element in the row space can be obtained by adding
some basis element, i.e. some row, to the latter. Then the problem of finding
a longest chain in the row space is to find as long a
sequence of edges as possible such that each new edge either adds 
a new vertex or a new color.  In general it is then optimal to either
add only one new vertex with existing colors or one new color with 
existing vertices.

The count of edges occurring in the sequence 
will give the size of the resulting chain in the row space.

Now consider a graph made up entirely of disjoint cycles whose
edges are colored various colors. Suppose, by taking multiple edges
of some colors, that in each cycle only one color is unique (occurs
only on one edge of that particular cycle, though it may occur
on other cycles), and
these unique colors are all different from each other.

We assert that optimum chains can
be made to have the following form: we add all but one edge in each
cycle and then we add the remaining cycle if it has a new color.
For, we can add all but one edge proceeding in sequence around the 
cycle, so that each edge adds a new vertex, 
and so the sum gives a new point in the row
space. And we can then add the final edge in that 
cycle with an increase if and only if it has a new color.
If we leave out more than one edge in a cycle, we might gain
some new colors for final edges elsewhere, but the net effect will be no better than
filling in all but one edge, since at most one edge is gained, one is lost
in so doing.

We may as well fill in the entire cycles in the order in which
we fill in their final edges, and then later do any cycles in
arbitrary order, for which we do not fill in their final cycles.
That is, the order of filling in non-final edges among different
cycles does not matter, and we may as well 
have as few as possible of them at each stage.
So then a cycle need not be started until the final edges prior to its 
final edge are put in place.

Now form a matrix to represent the possible priorities.  Among the cycles
which are to have some final edge, we can allow cycle $i$ to precede cycle $j$ if and
only if cycle $i$ does not contain the unique color for cycle $j$ (this color will be
multiple for cycle $i$, and off the main diagonal if it exists).  In this
case make the $(i,j)$ entry of the matrix equal to $1$.  This matrix can be arbitrary
except for $0$ on the main diagonal, by choosing the non-unique colors on each
cycle.

Then we can have final edges $k_1, k_2, \ldots k_n$ inserted if and only if
all the $(k_i,k_j)$ entries of this matrix are $1$, for $i<j$.  This means existence
of a triangular submatrix in the sense mentioned above (with some convention
on the main diagonal), and by the
lemma it is NP-complete. $\Box$

\section{GCD and LCM}

In any tropical division problem of dividing polynomials in one variable, 
$d$ into $s$, we can produce a least possible superquotient 
$q$ such that $dq$ is greater than or equal to $s$ just by the fact that the 
minimum of two superquotients is also a superquotient by distributivity
of products over minimum.

The least superquotient in dividing $\sum_0^n a_ix^i$ into $\sum_0^m b_ix_i$ can be computed 
explicitly by the requirements that its coefficients $c_i$ are the least numbers
which satisfy
$b_i\le \inf_k (a_{i-k}+c_k)$ so $c_k=\sup_i (b_i-a_{i-k})$.

Example. For two polynomials of the same degree, the least superquotient represents
the least translate raising one above the other.

For the same reason, any two nonnegative tropical polynomials 
will have a unique least nonnegative tropical common multiple in each degree
for which a common multiple exists, 
because the infimum of two common multiples is a common multiple, 
by the infimum of the corresponding quotients. A nonnegative 
tropical common multiple of finite degree is given by the product. 
The degree is chosen to be minimal.

In terms of degree there is a little paradox about minimal polynomials, 
insofar as a minimal polynomial of larger degree would be smaller, 
and so for this reason it may be preferable to fix the degrees in 
constructions involving superquotients.

We believe the following algorithm will converge rapidly to the tropical
least common multiple, or else it will be possible to tell rapidly that
convergence will not occur in the given degree, 
but we do not have a proof that it is polynomial
time in the degree and sizes of the coefficients. 
In fact over the rationals one can note denominators are bounded, and
it is likely that in projective space, the system remains within a bounded
region which must become eventually periodic.  In fact the maximum slope of the 
convex hull at each
stage should be bounded by that in $f$ and $g$: for products this follows
by Prop.1 and for superquotients by a similar argument.  Given a fixed range of 
degrees this means modulo translation that the construction stays within a fixed
rectangle.

Assume the constant terms are $0$.  We try different possibilities for the degree 
of the L.C.M.  In each case we then write out the initial trial $h$ for the L.C.M. 
as the tropical polynomial which has all coefficients $0$ up to this degree; this
will be a lower bound.  Repeat the 
following process until it converges or can be seen to increase
at a constant rate indefinitely: take the minimal superquotient of $h$ by $f$, and 
then multiply the result by $f$.  Take that as the new $h$; the L.C.M must be at least 
that great by monotonicity of superquotient.  Then take the minimal superquotient of 
the new $h$ by $g$, and multiply that by $g$ take that as the new superquotient. 
Coefficients never decrease in this process. At 
each step before convergence, at least one coefficient must increase.
If there is a common multiple in this degree, step by step the common multiple should
remain above this.

Definition. The tropical g.c.d. of two nonnegative tropical 
polynomials $f,g$ in one variable is the least superquotient 
of $fg$ by their least common multiple.

The g.c.d of polynomials over a tropical semiring is not in general unique.
We can take an example of nonunique factorization of 
some polynomial $f$ and compare $f$ with the product of the 
two factors.  For example, with Boolean polynomials, 
both $1+x^n$ and $1+x^{n+1}$  divide $f=1+x+x^2+\ldots +x^{2n+1}$ 
but their product uniquely factors and does not divide $f$.

\section{Algorithms for factoring}

In the Boolean case one reasonable choice of a factoring 
algorithm is a straightforward branch and bound algorithm 
based on choosing the two factors in degrees starting with $1$ and 
increasing, and stopping any branch when a product is not contained 
in the polynomial $g(x)$ to be factored.  As far as we know this
may have exponential time in the worst case; it would be of some interest
in this respect to know what Boolean prime polynomials of given degree
have the maximum number of ones.  
If $g(x)$ is either sparse or if the set of degrees missing from 
it is sparse, then this algorithm can probably be improved.  In 
the sparse case one might compare $g(x)$ with its translates to 
see when the intersections are large, and this gives a 
likely degree for some term in a factor.  In the case 
where the complementary set is sparse,   for each sum of a 
pair which could produce a missing degree, choose one of the two 
summands to be missing from a factor, as consistently as 
possible.  This is likely to produce a factorization, and 
can be made into a branch and bound algorithm with at most exponential time.

For factoring tropical polynomials, it seems from the above, 
that one will have to consider degrees that are not too large; perhaps 
degree $20$ might be the most that could reasonably be done on a 
personal computer in the general case.  If the degree is fixed, 
then the problem of factoring becomes a straightforward question 
of linear programming, which can be done in polynomial time using 
a homotopy algorithm over $Q^+$.  Integer linear programming is 
NP-complete, and it is not clear whether for fixed degree, factorization 
using only integer coefficients can be done in polynomial time, 
though the linear inequalities here are comparatively simple.

Example. This example is to show that for equal degree equal 
slope concave factorizations, it can happen that an integer 
polynomial is the square of a fractional polynomial but not an integer polynomial.

Consider first a case where $a_0=a_{200}=0,a_5=a_{105}=a_6=2,
a_{10}=5,a_9=4$, and whenever any coefficient except  $a_0,a_{10}$ 
is less than $24$ the coefficients in degree exactly $5$ higher is 
$4$ greater than that coefficient until a value of at least 
$20$ is reached, and all other coefficients are $24$.   
The square of this is a degree $400$ polynomial which is totally 
even. However in any square factorization the coefficient $c_{20} =10$ 
must be produced as the square of the $a_{10}$ term; $c_{30}$  and  
$c_{40}$ being $24$, as well as $c_j, j=1,2,3,4,7,8,211,212,213,
214$ imply there is no other product into degree $20$ 
which could produce a value of $10$ than  $a_5+b_5$. These calculations
can be checked by computer.

It seems however that if a factorization of a tropical polynomial over Z exists over Q, not requiring 
the factors to be equal, then there will be a factorization of the polynomial over Z.  This factorization 
is obtained simply by rounding all fractions up to the next higher integer in one factor, and rounding 
them down to the next lower integer in the other.  Suppose that a term in the old 
product is given by the sum of two fractional entries, then it is the
sum of the integer parts plus an integer which is a 
sum of two positive fractions which are less than 1, so the sum of the fractional
parts is 1, which the new 
procedure also gives for that sum.  Any other sum of one or two fractional entries 
which previously was at least this great will still be, because a sum of two positive 
fractions less than 1, is replaced by 1, and if that is a decrease, it will be the next 
lower integer.  Suppose a term in the old product is the sum of two integer entries in 
the old factorization, and is not a sum of two fractional entries.  The only way this rounding 
procedure could affect that is if that integer was less that the sum of the fractional 
entries previously, but now exceeds them.  But again a sum of two positive fractions 
less than 1 is replaced by 1; if this is a decrease it must be the next 
lower integer, and there can be no change.

We propose a branch and bound algorithm that in each degree $k$ 
of the product polynomial in turn chooses one of the possible sums 
$a_i+b_{k-i}$ to be $c_k$.  It uses this information to deduce 
inequalities about the $a,b$; and uses that to either terminate 
a branch or restrict future choices when possible.  It might be 
simplest to consider the concave case.  Then it will be best to start 
simultaneously from both extreme ends, since there will be few 
choices, and these will help determine later choices.  There will also be 
fewer choices for those coefficients $c$ which are closest above the 
convex hull of the diagram.  
This algorithm will be conceivably factorial time in the degree, and 
polynomial in the size of the coefficients given that degree.

The following result is perhaps not what it seems, but can be combined 
with other algorithms to improve them.

Theorem 14.  For tropical polynomials of any given positive degree $n$ in 
each of $x,y$ separately, with all coefficients up to the maximal degree in each variable
being present (finite) there is an algorithm which provably produces an optimal 
factorization with probability $1$ in polynomial time, that is, the exceptional 
set is a finite union of algebraic sets of codimension at least $1$, and a 
number of arithmetic operations is involved which is polynomial in $n$.

Proof.  It suffices to give a polynomial time algorithm which will prove 
irreducibility except on such a codimension $1$ subset. But this can be done
trivially by looking at the 4 corners, degrees $(0,0),(0,2n),(2n,0),(2n,2n).\Box$

A study of the convex case and convex hulls in general should produce better
criteria.  In particular, it seems that if one flattens out the convex hulls
in particular directions, the flattened and simplified convex hulls must also
factor.
If one has total concavity then the 8 outer edges of the product determine the 
4 edges of the factors up to choice of which is in which factor, and this should 
determine inner edges.

\section{Conclusion}

In the general case factorization of polynomials even in one variable over 
the tropical algebra is quite difficult and algorithms will be practical 
only for small degrees or in special cases.  This is true even if 
coefficients are $0$ and infinity.  For bounded degree, factorization 
is quite possible, and many irreducibility conditions can be derived.
Problems of both factorization and tropical rank are NP-complete.

There is a polynomial time parallel algorithm for factorization 
over the tropical semiring, since parallel polynomial time is equivalent to 
PSPACE which contains NP.  (Any factorization must wind up with 
factors which as numbers, involve a polynomial number of digits).  
If, as conjectured, quantum P time is strictly a subset of NP there 
cannot be a quantum polynomial time algorithm for factoring tropical
or Boolean polynomials, unlike for the problem 
of factoring integers.    

We have given algorithms which in a 
certain sense statistically factor most tropical polynomials 
in $2$ variables polynomial time, by proving them irreducible.  We have
suggested directions for branch and bound algorithms to extend these.

For Boolean two variable polynomials $h$ which are products of polynomials 
$f,g$ each having degree $n$ in $x,y$ separately, a simple counting argument, 
unlike in the $1$ variable case, shows that the proportion of irreducible 
polynomials tends to $1$.  It would be of interest to find a simple criterion which 
in polynomial time can prove irreducibility with probability approaching $1$ 
as $n\rightarrow \infty$.

The average polynomial $h$ will have about half of its coefficients 1 and the
other half zero.  For these, it will probably be the case that $f,g$ must
have a much smaller proportion of their coefficients 1, something like the
square roots of their degrees.  It might be possible to prove some result
of this kind using existing combinatorial techniques, but it seems far
from trivial.  Given this, knowledge of the sparsity of nonzero terms in $f,g$,
a branch and bound algorithm which starts with the lowest degrees and assigns
coefficients in $h$ (assumed to have constant term 1) alternately to $f$
and to $g$ might proceed somewhat rapidly because of the high likelihood of finding
a product term which does not actually exist in $h$.  If $B_k$ branches survive to
level $k$ (which is the number of terms found so far) 
for $h$ of degree $n$ in $x,y$, then the number of branches which survive
to the next level can be estimated as $B_k 2^{-[k/2]}(n^2/2-k+1)$ since we can choose the
new term in about $n^2/2-k+1$ ways, allowing that the current
term in $f$ could be repeated in $g$ and $[k/2]$ new products with it must be a 1 coefficient,
and the probability of this happening could be estimated as $1/2$.  This
should take time about $n^{C\log n}$ to complete since as soon as $2^k$ exceeds
$(n^2)^2$ the series will converge rapidly.  

However this algorithm and estimate can be greatly improved by choosing the terms in order
according to the maximum of the degrees in $x,y$ separately, then the minimum
degree and variable, and moreover using the previous 
term in the factor as a lower bound, and the
next unfactored term in the product as an upper bound.  In this case, as soon as
we have about degree $C(\log (n))^2$ in the product we would expect to have
$C_1 \log (n)$ terms in the factors, and we might expect $k$ could be
taken about $\log \log (n)$. This algorithm could be polynomial time in
most cases. 
It would be however difficult to prove this behavior without making
some statistical hypotheses.

Open problems:

1. Is a random Boolean polynomial of large degree in 1 variable more likely to be
prime or factorable?

2. Is there a polynomial time algorithm to find a nontrivial common factor of two 
tropical polynomials, given that such a factor exists? Both eliminants and
L.C.M. give necessary conditions.

Those readers interested in other work on tropical algebra should see the
website www.arxiv.org where this paper was posted; additional open problems
are listed at www.math.umn.edu/$\sim$develin/tropicalproblems.html.

\end{document}